\documentclass[10pt]{article}
\usepackage{cite}
\usepackage{mathrsfs}
\usepackage{amsfonts}
\usepackage{amsmath}
\usepackage{amsfonts,amssymb,color}
\usepackage{dsfont}
\usepackage{curves}
\usepackage{mathrsfs}
\usepackage{pifont}
\usepackage{amssymb}
\allowdisplaybreaks

\numberwithin{equation}{section}

\date{}

\def\BigRoman{\uppercase\expandafter{\romannumeral\number\count 255 }}
\def\Romannumeral{\afterassignment\BigRoman\count255=}

\begin{document}
\title{Sufficient conditions for a graph with minimum degree to have a component factor
%\thanks{}
}
\author{\small  Jie Wu\footnote{Corresponding
author. E-mail address: jskjdxwj@126.com}\\
\small  School of Economics and management,\\
\small  Jiangsu University of Science and Technology,\\
\small  Zhenjiang, Jiangsu 212100, China\\
}

\maketitle
\begin{abstract}
\noindent Let $\mathcal{T}_{\frac{k}{r}}$ denote the set of trees $T$ such that $i(T-S)\leq\frac{k}{r}|S|$ for any $S\subset V(T)$ and for
any $e\in E(T)$ there exists a set $S^{*}\subset V(T)$ with $i((T-e)-S^{*})>\frac{k}{r}|S^{*}|$, where $r<k$ are two positive integers. A
$\{C_{2i+1},T:1\leq i<\frac{r}{k-r},T\in\mathcal{T}_{\frac{k}{r}}\}$-factor of a graph $G$ is a spanning subgraph of $G$, in which every
component is isomorphic to an element in $\{C_{2i+1},T:1\leq i<\frac{r}{k-r},T\in\mathcal{T}_{\frac{k}{r}}\}$. Let $A(G)$ and $Q(G)$ denote
the adjacency matrix and the signless Laplacian matrix of $G$, respectively. The adjacency spectral radius and the signless Laplacian spectral
radius of $G$, denoted by $\rho(G)$ and $q(G)$, are the largest eigenvalues of $A(G)$ and $Q(G)$, respectively. In this paper, we study the
connections between the spectral radius and the existence of a $\{C_{2i+1},T:1\leq i<\frac{r}{k-r},T\in\mathcal{T}_{\frac{k}{r}}\}$-factor
in a graph. We first establish a tight sufficient condition involving the adjacency spectral radius to guarantee the existence of a
$\{C_{2i+1},T:1\leq i<\frac{r}{k-r},T\in\mathcal{T}_{\frac{k}{r}}\}$-factor in a graph. Then we propose a tight signless Laplacian spectral
radius condition for the existence of a $\{C_{2i+1},T:1\leq i<\frac{r}{k-r},T\in\mathcal{T}_{\frac{k}{r}}\}$-factor in a graph.
\\
\begin{flushleft}

(2020) Mathematics Subject Classification: 05C50, 05C70, 90B99
\end{flushleft}
\end{abstract}

\section{Introduction}

In this paper, we deal with finite and undirected graphs which have neither loops nor multiple edges. Let $G$ be a graph. We denote by $V(G)$
and $E(G)$ the set of vertices and the set of edges of $G$, respectively. The order of $G$ is the number $n=|V(G)|$ of its vertices. The size
of $G$ is the number $e(G)=|E(G)|$ of its edges. For $v\in V(G)$, the degree of $v$ in $G$ is denoted by $d_G(v)$. Let $i(G)$ and $\delta(G)$
denote the number of isolated vertices and the minimum degree of $G$, respectively. For any $S\subseteq V(G)$, $G[S]$ is the subgraph of $G$
induced by $S$ and $G-S$ is the subgraph of $G$ induced by $V(G)-S$. For any $E'\subseteq E(G)$, let $G-E'$ denote the subgraph obtained from
$G$ by deleting $E'$. For convenience, write $G-v=G-\{v\}$ for $S=\{v\}$ and $G-e=G-\{e\}$ for $E'=\{e\}$. As usual, let $P_n$, $C_n$, $K_n$
and $K_{1,n-1}$ denote the path, the circuit, the complete graph and the star of order $n$, respectively. A connected graph without circuits
is called a tree, which is denoted by $T$. For any two positive integers $k$ and $r$ with $r<k$, let $\mathcal{T}_{\frac{k}{r}}$ denote the
set of trees $T$ such that $i(T-S)\leq\frac{k}{r}|S|$ for any $S\subset V(T)$ and for any $e\in E(T)$ there exists a set $S^{*}\subset V(T)$
with $i((T-e)-S^{*})>\frac{k}{r}|S^{*}|$. For two given graphs $G_1$ and $G_2$, we denote by $G_1\cup G_2$ their union. The join $G_1\vee G_2$
is obtained from $G_1\cup G_2$ by joining each vertex of $G_1$ with each vertex of $G_2$ by an edge. Let $c$ be a real number. Recall that
$\lfloor c\rfloor$ is the greatest integer with $\lfloor c\rfloor\leq c$.

Let $\mathcal{H}$ denote a set of connected graphs. A subgraph $H$ of $G$ is called an $\mathcal{H}$-factor of $G$ if $V(H)=V(G)$ and each
component of $H$ is isomorphic to an element of $\mathcal{H}$. An $\mathcal{H}$-factor is also referred as a component factor. An
$\mathcal{H}$-factor is called a $P_{\geq k}$-factor if $\mathcal{H}=\{P_k,P_{k+1},\ldots\}$. An $\mathcal{H}$-factor is called a
$\{C_{2i+1},T:1\leq i<\frac{r}{k-r},T\in\mathcal{T}_{\frac{k}{r}}\}$-factor if
$\mathcal{H}=\{C_{2i+1},T:1\leq i<\frac{r}{k-r},T\in\mathcal{T}_{\frac{k}{r}}\}$. An $\mathcal{H}$-factor means a star-factor in which every
component is a star. Note that a perfect matching is indeed a $\{P_2\}$-factor of $G$.

Kaneko \cite{Ka} established a criterion for a graph with a $P_{\geq3}$-factor. Liu and Pan \cite{LP}, Dai \cite{Dai}, Wu \cite{Wp} provided
some sufficient conditions for the existence of $P_{\geq3}$-factors in graphs. Ando et al \cite{AEKKM} proved that a claw-free graph with
minimum degree at least $d$ contains a $P_{\geq d+1}$-factor. Tutte \cite{T} showed a necessary and sufficient condition for a graph to have
a $\{K_2,C_i:i\geq3\}$-factor. Klopp and Steffen \cite{KSf} investigated the existence of $\{K_{1,1},K_{1,2},C_i:i\geq3\}$-factors in graphs.
Zhou, Xu and Sun \cite{ZXS} proposed some sufficient conditions for graphs to contain $\{K_{1,j}:1\leq j\leq k\}$-factors. Kano and Saito
\cite{KS} verified that a graph $G$ satisfying $i(G-S)\leq\frac{1}{k}|S|$ for any $S\subset V(G)$ has a $\{K_{1,j}:k\leq j\leq2k\}$-factor.
Kano, Lu and Yu \cite{KLY} provided sufficient conditions using isolated vertices for component factors with every component of order at least
three and proved that a graph $G$ satisfying $i(G-S)\leq\frac{|S|}{2}$ for any $S\subset V(G)$ contains a $\{K_{1,2},K_{1,3},K_5\}$-factor.
Wolf \cite{Wolf} claimed a characterization using isolated vertices for a graph with a
$\{C_{2i+1},T:1\leq i<\frac{r}{k-r},T\in\mathcal{T}_{\frac{k}{r}}\}$-factor. For other sufficient conditions for the existence of graph factors
in graphs, see \cite{GWC,ZPX1,Zr,Zt,ZZL,Wr,Wa}.

Given a graph $G$ with vertex set $V(G)=\{v_1,v_2,\ldots,v_n\}$, the adjacency matrix $A(G)=(a_{ij})_{n\times n}$ of $G$ is a 0--1 matrix
in which the entry $a_{ij}=1$ if and only if $v_iv_j\in E(G)$. Let $D(G)$ denote the diagonal matrix of vertex degrees of $G$. The signless
Laplacian matrix $Q(G)$ of $G$ are defined by $Q(G)=D(G)+A(G)$. The largest eigenvalue of $A(G)$ is called the adjacency spectral radius of
$G$, denoted by $\rho(G)$. The largest eigenvalue of $Q(G)$ is called the signless Laplacian spectral radius of $G$, denoted by $q(G)$.

O \cite{O}, Zhou, Sun and Zhang \cite{ZSZ} proved two sharp upper bounds for the adjacency spectral radius in a graph without a $\{P_2\}$-factor.
Zhou and Zhang \cite{ZZ} gave a lower bound on the signless Laplacian spectral radius of $G$ to guarantee that $G$ contains a $\{P_2\}$-factor.
Zhou, Sun and Liu \cite{ZSL2}, Zhou, Zhang and Sun \cite{ZZS} presented two spectral radius conditions for graphs to possess $P_{\geq2}$-factors.
Wu \cite{Wc}, Wang and Zhang \cite{WZ}, Zhou, Sun and Liu \cite{ZSL1}, Zhou and Wu \cite{ZW} provided some spectral radius conditions for the
existence of spanning trees in connected graphs. Zhou \cite{Zs} proposed two spectral radius conditions for bipartite graphs to have star-factors.
Zhou and Liu \cite{ZL1} put forward a lower bound on the $A_a$-spectral radius for a connected graph to possess a $\{K_{1,j}:m\leq j\leq2m\}$-factor.
Lv, Li and Xu \cite{LLX} showed a sufficient condition involving the $A_{\alpha}$-spectral radius for a graph to have a $\{K_2,C_{2i+1}: i\geq1\}$-factor,
and gave a distance signless Laplacian spectral radius condition for a graph to have a $\{K_2,C_{2i+1}: i\geq1\}$-factor. Miao and Li \cite{ML}
obtained some sufficient conditions involving the adjacency spectral radius and the distance spectral radius for the existence of
$\{K_{1,j}:1\leq j\leq k\}$-factors in graphs.

Motivated by \cite{Wolf} directly, we first propose an adjacency spectral radius condition for a connected graph with minimum degree $\delta$ to have
a $\{C_{2i+1},T:1\leq i<\frac{r}{k-r},T\in\mathcal{T}_{\frac{k}{r}}\}$-factor, then we obtain a signless Laplacian spectral radius condition for a
connected graph with minimum degree $\delta$ to have a $\{C_{2i+1},T:1\leq i<\frac{r}{k-r},T\in\mathcal{T}_{\frac{k}{r}}\}$-factor.

\medskip

\noindent{\textbf{Theorem 1.1.}} Let $k$ and $r$ be two positive integers with $r<k$, and let $G$ be a connected graph of order $n$ with $\delta(G)=\delta$
and $n\geq\max\Big\{\frac{(k+r)(k+2r)(k\delta+k+r)}{k^{2}r},\frac{2kr\delta^{2}+(2k^{2}+kr+2r^{2})\delta+k^{2}+3kr-2r^{2}}{2r(k-r)}\Big\}$. If
$$
\rho(G)\geq\rho\Big(K_{\delta}\vee\Big(K_{n-\lfloor\frac{k\delta}{r}\rfloor-\delta-1}\cup\Big(\Big\lfloor\frac{k\delta}{r}\Big\rfloor+1\Big)K_1\Big)\Big),
$$
then $G$ has a $\{C_{2i+1},T:1\leq i<\frac{r}{k-r},T\in\mathcal{T}_{\frac{k}{r}}\}$-factor unless
$G=K_{\delta}\vee(K_{n-\lfloor\frac{k\delta}{r}\rfloor-\delta-1}\cup(\lfloor\frac{k\delta}{r}\rfloor+1)K_1)$.

\medskip

\noindent{\textbf{Theorem 1.2.}} Let $k$ and $r$ be two positive integers with $r<k$, and let $G$ be a connected graph of order $n$ with $\delta(G)=\delta$
and $n\geq\max\Big\{\frac{(k+r)(k+2r)(k\delta+k+r)}{k^{2}r},\frac{(k^{2}+2kr)\delta^{2}+(2k^{2}+3kr+2r^{2})\delta+k^{2}+3kr}{2r(k-r)}\Big\}$. If
$$
q(G)\geq q\Big(K_{\delta}\vee\Big(K_{n-\lfloor\frac{k\delta}{r}\rfloor-\delta-1}\cup\Big(\Big\lfloor\frac{k\delta}{r}\Big\rfloor+1\Big)K_1\Big)\Big),
$$
then $G$ has a $\{C_{2i+1},T:1\leq i<\frac{r}{k-r},T\in\mathcal{T}_{\frac{k}{r}}\}$-factor unless
$G=K_{\delta}\vee(K_{n-\lfloor\frac{k\delta}{r}\rfloor-\delta-1}\cup(\lfloor\frac{k\delta}{r}\rfloor+1)K_1)$.

\medskip

\section{Preliminary lemmas}

In this section, we show some lemmas, which will be used to verify our main results. Wolf \cite{Wolf} claimed a characterization for a graph with
a $\{C_{2i+1},T:1\leq i<\frac{r}{k-r},T\in\mathcal{T}_{\frac{k}{r}}\}$-factor.

\medskip

\noindent{\textbf{Lemma 2.1}} (Wolf \cite{Wolf}). Let $k$ and $r$ be two positive integers with $r<k$, and let $G$ be a graph. Then $G$
has a $\{C_{2i+1},T:1\leq i<\frac{r}{k-r},T\in\mathcal{T}_{\frac{k}{r}}\}$-factor if and only if
$$
i(G-S)\leq\frac{k}{r}|S|
$$
for any $S\subset V(G)$.

\medskip

\noindent{\textbf{Lemma 2.2}} (Li and Feng \cite{LF}). Let $G$ be a connected graph and let $H$ be a subgraph of $G$. Then
$$
\rho(G)\geq\rho(H),
$$
with equality if and only if $G=H$.

\medskip

\noindent{\textbf{Lemma 2.3}} (Hong \cite{Ha}). Let $G$ be a graph with $n$ vertices. Then
$$
\rho(G)\leq\sqrt{2e(G)-n+1},
$$
where the equality holds if and only if $G$ is a star or a complete graph.

\medskip

\noindent{\textbf{Lemma 2.4}} (Shen, You, Zhang and Li \cite{SYZL}). Let $G$ be a connected graph. If $H$ is a subgraph of $G$, then
$$
q(G)\geq q(H),
$$
with equality holding if and only if $G=H$.

\medskip

\noindent{\textbf{Lemma 2.5}} (Das \cite{Dm}). Let $G$ be a graph of order $n$. Then
$$
q(G)\leq\frac{2e(G)}{n-1}+n-2.
$$

\medskip

\section{The proof of Theorem 1.1}

\noindent{\it Proof of Theorem 1.1.} Assume that $G$ has no $\{C_{2i+1},T:1\leq i<\frac{r}{k-r},T\in\mathcal{T}_{\frac{k}{r}}\}$-factor. By virtue
of Lemma 2.1, there exists some nonempty subset $S$ of $V(G)$ such that
$$
i(G-S)>\frac{k}{r}|S|.
$$
In terms of the integrity of $i(G-S)$, we possess
$$
i(G-S)\geq\Big\lfloor\frac{k}{r}|S|\Big\rfloor+1.
$$
Let $|S|=s$. Then $G$ is a spanning subgraph of $G_1=K_s\vee(K_{n-\lfloor\frac{ks}{r}\rfloor-s-1}\cup(\lfloor\frac{ks}{r}\rfloor+1)K_1)$. Together
with Lemma 2.2, we deduce
\begin{align}\label{eq:3.1}
\rho(G)\leq\rho(G_1),
\end{align}
where the equality holds if and only if $G=G_1$. Notice that $\delta(G)=\delta$ and $\delta(G_1)\geq\delta(G)$. Thus, we get
$s=\delta(G_1)\geq\delta(G)=\delta$. The following proof will be divided into two cases according to the value of $s$.

\medskip

\noindent{\bf Case 1.} $s=\delta$.

In this case, $G_1=K_{\delta}\vee(K_{n-\lfloor\frac{k\delta}{r}\rfloor-\delta-1}\cup(\lfloor\frac{k\delta}{r}\rfloor+1)K_1)$. Together with \eqref{eq:3.1},
we conclude
$$
\rho(G)\leq\rho\Big(K_{\delta}\vee\Big(K_{n-\lfloor\frac{k\delta}{r}\rfloor-\delta-1}\cup\Big(\Big\lfloor\frac{k\delta}{r}\Big\rfloor+1\Big)K_1\Big)\Big),
$$
with equality holding if and only if $G=K_{\delta}\vee(K_{n-\lfloor\frac{k\delta}{r}\rfloor-\delta-1}\cup(\lfloor\frac{k\delta}{r}\rfloor+1)K_1)$. Observe
that $K_{\delta}\vee(K_{n-\lfloor\frac{k\delta}{r}\rfloor-\delta-1}\cup(\lfloor\frac{k\delta}{r}\rfloor+1)K_1)$ has no $\{C_{2i+1},T:1\leq i<\frac{r}{k-r},T\in\mathcal{T}_{\frac{k}{r}}\}$-factor. Thus, we can get a contradiction.

\medskip

\noindent{\bf Case 2.} $s\geq\delta+1$.

Recall that $G_1=K_s\vee(K_{n-\lfloor\frac{ks}{r}\rfloor-s-1}\cup(\lfloor\frac{ks}{r}\rfloor+1)K_1)$. By virtue of Lemma 2.3,
$\frac{ks}{r}-1<\lfloor\frac{ks}{r}\rfloor\leq\frac{ks}{r}$ and $n\geq\lfloor\frac{ks}{r}\rfloor+s+1>\frac{ks}{r}+s$, we obtain
\begin{align}\label{eq:3.2}
\rho(G_1)\leq&\sqrt{2e(G_1)-n+1}\nonumber\\
=&\sqrt{2\binom{n-\lfloor\frac{ks}{r}\rfloor-1}{2}+2s\Big(\Big\lfloor\frac{ks}{r}\Big\rfloor+1\Big)-n+1}\nonumber\\
=&\sqrt{\Big(n-\Big\lfloor\frac{ks}{r}\Big\rfloor-1\Big)\Big(n-\Big\lfloor\frac{ks}{r}\Big\rfloor-2\Big)+2s\Big(\Big\lfloor\frac{ks}{r}\Big\rfloor+1\Big)-n+1}\nonumber\\
<&\sqrt{\Big(n-\Big(\frac{ks}{r}-1\Big)-1\Big)\Big(n-\Big(\frac{ks}{r}-1\Big)-2\Big)+2s\Big(\Big\lfloor\frac{ks}{r}\Big\rfloor+1\Big)-n+1}\nonumber\\
=&\frac{1}{r}\sqrt{(k^{2}+2kr)s^{2}-(2krn-2r^{2}-kr)s+r^{2}n^{2}-2r^{2}n+r^{2}}.
\end{align}
Let $f(s)=(k^{2}+2kr)s^{2}-(2krn-2r^{2}-kr)s+r^{2}n^{2}-2r^{2}n+r^{2}$. Since $n\geq\lfloor\frac{ks}{r}\rfloor+s+1>\frac{ks}{r}+s$, we possess
$\delta+1\leq s<\frac{rn}{k+r}$. By a direct computation, we get
\begin{align*}
f(\delta+1)-f\Big(\frac{rn}{k+r}\Big)=&(k^{2}+2kr)(\delta+1)^{2}-(2krn-2r^{2}-kr)(\delta+1)+r^{2}n^{2}-2r^{2}n+r^{2}\\
&-\Big((k^{2}+2kr)\Big(\frac{rn}{k+r}\Big)^{2}-(2krn-2r^{2}-kr)\Big(\frac{rn}{k+r}\Big)+r^{2}n^{2}-2r^{2}n+r^{2}\Big)\\
=&\Big(\frac{rn}{k+r}-\delta-1\Big)\Big(\frac{k^{2}rn}{k+r}-(k+2r)(k\delta+k+r)\Big)\\
>&0,
\end{align*}
where the inequality holds from the fact that
\begin{align*}
n>&\max\Big\{\frac{(k+r)(k+2r)(k\delta+k+r)}{k^{2}r},\frac{2kr\delta^{2}+(2k^{2}+kr+2r^{2})\delta+k^{2}+3kr-2r^{2}}{2r(k-r)}\Big\}\\
\geq&\frac{(k+r)(k+2r)(k\delta+k+r)}{k^{2}r}\\
>&\frac{(k+r)(\delta+1)}{r}.
\end{align*}
This implies that, for $\delta+1\leq s<\frac{rn}{k+r}$, the function $f(s)$ attains its maximum value at $s=\delta+1$. Combining this with
\eqref{eq:3.1}, \eqref{eq:3.2} and
$n>\max\Big\{\frac{(k+r)(k+2r)(k\delta+k+r)}{k^{2}r},\frac{2kr\delta^{2}+(2k^{2}+kr+2r^{2})\delta+k^{2}+3kr-2r^{2}}{2r(k-r)}\Big\}
\geq\frac{2kr\delta^{2}+(2k^{2}+kr+2r^{2})\delta+k^{2}+3kr-2r^{2}}{2r(k-r)}$, we obtain
\begin{align}\label{eq:3.3}
\rho(G)\leq&\rho(G_1)\nonumber\\
<&\frac{1}{r}\sqrt{f(\delta+1)}\nonumber\\
=&\frac{1}{r}\sqrt{(k^{2}+2kr)(\delta+1)^{2}-(2krn-2r^{2}-kr)(\delta+1)+r^{2}n^{2}-2r^{2}n+r^{2}}\nonumber\\
=&\frac{1}{r}\sqrt{(rn-k\delta-2r)^{2}-2r(k-r)n+2kr\delta^{2}+(2k^{2}+kr+2r^{2})\delta+k^{2}+3kr-2r^{2}}\nonumber\\
<&\frac{1}{r}(rn-k\delta-2r).
\end{align}

Since $K_{n-\lfloor\frac{k\delta}{r}\rfloor-1}$ is a proper subgraph of
$K_{\delta}\vee(K_{n-\lfloor\frac{k\delta}{r}\rfloor-\delta-1}\cup(\lfloor\frac{k\delta}{r}\rfloor+1)K_1)$, it follows from Lemma 2.2,
$\lfloor\frac{k\delta}{r}\rfloor\leq\frac{k\delta}{r}$ and the hypothesis of the theorem that
\begin{align*}
\rho(G)\geq&\rho\Big(K_{\delta}\vee\Big(K_{n-\lfloor\frac{k\delta}{r}\rfloor-\delta-1}\cup\Big(\Big\lfloor\frac{k\delta}{r}\Big\rfloor+1\Big)K_1\Big)\Big)\\
>&\rho(K_{n-\lfloor\frac{k\delta}{r}\rfloor-1})\\
=&n-\Big\lfloor\frac{k\delta}{r}\Big\rfloor-2\\
\geq&n-\frac{k\delta}{r}-2\\
=&\frac{1}{r}(rn-k\delta-2r),
\end{align*}
which leads to a contradiction to \eqref{eq:3.3}. Theorem 1.1 is proved. \hfill $\Box$

\section{The proof of Theorem 1.2}

\noindent{\it Proof of Theorem 1.2.} Assume that $G$ has no $\{C_{2i+1},T:1\leq i<\frac{r}{k-r},T\in\mathcal{T}_{\frac{k}{r}}\}$-factor. Then using
Lemma 2.1, there exists some nonempty subset $S$ of $V(G)$ such that
$$
i(G-S)>\frac{k}{r}|S|.
$$
According to the integrity of $i(G-S)$, we obtain
$$
i(G-S)\geq\Big\lfloor\frac{k}{r}|S|\Big\rfloor+1.
$$
Let $|S|=s$. Then $G$ is a spanning subgraph of $G_1=K_s\vee(K_{n-\lfloor\frac{ks}{r}\rfloor-s-1}\cup(\lfloor\frac{ks}{r}\rfloor+1)K_1)$. Together
with Lemma 2.4, we possess
\begin{align}\label{eq:4.1}
q(G)\leq q(G_1),
\end{align}
where the equality holds if and only if $G=G_1$. Note that $\delta(G)=\delta$ and $\delta(G_1)=s\geq\delta(G)$. Thus, we get $s\geq\delta$. In what
follows, we shall consider two cases by the value of $s$.

\medskip

\noindent{\bf Case 1.} $s=\delta$.

In this case, $G_1=K_{\delta}\vee(K_{n-\lfloor\frac{k\delta}{r}\rfloor-\delta-1}\cup(\lfloor\frac{k\delta}{r}\rfloor+1)K_1)$. In terms of \eqref{eq:4.1},
we obtain
$$
q(G)\leq q\Big(K_{\delta}\vee\Big(K_{n-\lfloor\frac{k\delta}{r}\rfloor-\delta-1}\cup\Big(\Big\lfloor\frac{k\delta}{r}\Big\rfloor+1\Big)K_1\Big)\Big),
$$
where the equality holds if and only if $G=K_{\delta}\vee(K_{n-\lfloor\frac{k\delta}{r}\rfloor-\delta-1}\cup(\lfloor\frac{k\delta}{r}\rfloor+1)K_1)$.
Observe that $K_{\delta}\vee(K_{n-\lfloor\frac{k\delta}{r}\rfloor-\delta-1}\cup(\lfloor\frac{k\delta}{r}\rfloor+1)K_1)$ contains no $\{C_{2i+1},T:1\leq i<\frac{r}{k-r},T\in\mathcal{T}_{\frac{k}{r}}\}$-factor. Thus, we can obtain a contradiction.

\medskip

\noindent{\bf Case 2.} $s\geq\delta+1$.

Recall that $G_1=K_s\vee(K_{n-\lfloor\frac{ks}{r}\rfloor-s-1}\cup(\lfloor\frac{ks}{r}\rfloor+1)K_1)$. It follows from Lemma 2.5,
$\frac{ks}{r}-1<\lfloor\frac{ks}{r}\rfloor\leq\frac{ks}{r}$ and $n\geq\lfloor\frac{ks}{r}\rfloor+s+1>\frac{ks}{r}+s$ that
\begin{align}\label{eq:4.2}
q(G_1)\leq&\frac{2e(G_1)}{n-1}+n-2\nonumber\\
=&\frac{2\binom{n-\lfloor\frac{ks}{r}\rfloor-1}{2}+2s\Big(\Big\lfloor\frac{ks}{r}\Big\rfloor+1\Big)}{n-1}+n-2\nonumber\\
=&\frac{\Big(n-\Big\lfloor\frac{ks}{r}\Big\rfloor-1\Big)\Big(n-\Big\lfloor\frac{ks}{r}\Big\rfloor-2\Big)+2s\Big(\Big\lfloor\frac{ks}{r}\Big\rfloor+1\Big)}{n-1}
+n-2\nonumber\\
<&\frac{\Big(n-\Big(\frac{ks}{r}-1\Big)-1\Big)\Big(n-\Big(\frac{ks}{r}-1\Big)-2\Big)+2s\Big(\frac{ks}{r}+1\Big)}{n-1}+n-2\nonumber\\
=&\frac{\Big(n-\frac{ks}{r}\Big)\Big(n-\frac{ks}{r}-1\Big)+2s\Big(\frac{ks}{r}+1\Big)}{n-1}+n-2\nonumber\\
=&\frac{(k^{2}+2kr)s^{2}-(2krn-kr-2r^{2})s+2r^{2}n^{2}-4r^{2}n+2r^{2}}{r^{2}(n-1)}.
\end{align}
Let $g(s)=(k^{2}+2kr)s^{2}-(2krn-kr-2r^{2})s+2r^{2}n^{2}-4r^{2}n+2r^{2}$. Since $n\geq\lfloor\frac{ks}{r}\rfloor+s+1>\frac{ks}{r}+s$, we deduce
$\delta+1\leq s<\frac{rn}{k+r}$. By a simple computation, we obtain
\begin{align*}
g(\delta+1)-g\Big(\frac{rn}{k+r}\Big)=&(k^{2}+2kr)(\delta+1)^{2}-(2krn-kr-2r^{2})(\delta+1)+2r^{2}n^{2}-4r^{2}n+2r^{2}\\
&-\Big((k^{2}+2kr)\Big(\frac{rn}{k+r}\Big)^{2}-(2krn-kr-2r^{2})\Big(\frac{rn}{k+r}\Big)+2r^{2}n^{2}-4r^{2}n+2r^{2}\Big)\\
=&\Big(\frac{rn}{k+r}-\delta-1\Big)\Big(\frac{k^{2}rn}{k+r}-(k+2r)(k\delta+k+r)\Big)\\
>&0,
\end{align*}
where the inequality holds from the fact that
\begin{align*}
n>&\max\Big\{\frac{(k+r)(k+2r)(k\delta+k+r)}{k^{2}r},\frac{(k^{2}+2kr)\delta^{2}+(2k^{2}+3kr+2r^{2})\delta+k^{2}+3kr}{2r(k-r)}\Big\}\\
\geq&\frac{(k+r)(k+2r)(k\delta+k+r)}{k^{2}r}\\
>&\frac{(k+r)(\delta+1)}{r}.
\end{align*}
This implies that, for $\delta+1\leq s<\frac{rn}{k+r}$, the function $g(s)$ attains its maximum value at $s=\delta+1$. Combining this with
\eqref{eq:4.1}, \eqref{eq:4.2} and
$n>\max\Big\{\frac{(k+r)(k+2r)(k\delta+k+r)}{k^{2}r},\frac{(k^{2}+2kr)\delta^{2}+(2k^{2}+3kr+2r^{2})\delta+k^{2}+3kr}{2r(k-r)}\Big\}
\geq\frac{(k^{2}+2kr)\delta^{2}+(2k^{2}+3kr+2r^{2})\delta+k^{2}+3kr}{2r(k-r)}$, we conclude
\begin{align}\label{eq:4.3}
q(G)\leq&q(G_1)\nonumber\\
<&\frac{g(\delta+1)}{r^{2}(n-1)}\nonumber\\
=&\frac{(k^{2}+2kr)(\delta+1)^{2}-(2krn-kr-2r^{2})(\delta+1)+2r^{2}n^{2}-4r^{2}n+2r^{2}}{r^{2}(n-1)}\nonumber\\
=&\frac{2(rn-k\delta-2r)}{r}-\frac{2r(k-r)n-(k^{2}+2kr)\delta^{2}-(2k^{2}+3kr+2r^{2})\delta-k^{2}-3kr}{r^{2}(n-1)}\nonumber\\
<&\frac{2(rn-k\delta-2r)}{r}.
\end{align}

Note that $K_{\delta}\vee(K_{n-\lfloor\frac{k\delta}{r}\rfloor-\delta-1}\cup(\lfloor\frac{k\delta}{r}\rfloor+1)K_1)$ contains
$K_{n-\lfloor\frac{k\delta}{r}\rfloor-1}$ as a proper subgraph. Together with Lemma 2.4, $\lfloor\frac{k\delta}{r}\rfloor\leq\frac{k\delta}{r}$
and the assumption of the theorem, we possess
\begin{align*}
q(G)\geq&q\Big(K_{\delta}\vee\Big(K_{n-\lfloor\frac{k\delta}{r}\rfloor-\delta-1}\cup\Big(\Big\lfloor\frac{k\delta}{r}\Big\rfloor+1\Big)K_1\Big)\Big)\\
>&q(K_{n-\lfloor\frac{k\delta}{r}\rfloor-1})\\
=&2\Big(n-\Big\lfloor\frac{k\delta}{r}\Big\rfloor-2\Big)\\
\geq&2\Big(n-\frac{k\delta}{r}-2\Big)\\
=&\frac{2(rn-k\delta-2r)}{r},
\end{align*}
which is to a contradiction to \eqref{eq:4.3}. This completes the proof of Theorem 1.2. \hfill $\Box$

\section{Concluding remarks}

In this paper, we provide two sufficient conditions to ensure that a connected graph $G$ has a
$\{C_{2i+1},T:1\leq i<\frac{r}{k-r},T\in\mathcal{T}_{\frac{k}{r}}\}$-factor in terms of its adjacency spectral radius and signless Laplacian
spectral radius. It is natural and interesting to propose some other spectral sufficient conditions to guarantee that a connected graph $G$ has
a $\{C_{2i+1},T:1\leq i<\frac{r}{k-r},T\in\mathcal{T}_{\frac{k}{r}}\}$-factor. It is also natural and interesting to put forward some spectral
sufficient conditions to ensure that a connected graph $G$ has some other substructure.

\section*{Declaration of competing interest}

\medskip

The author declares that he has no known competing financial interests or personal relationships that could have appeared to influence the work
reported in this paper.

\section*{Data availability}

\medskip

No data was used for the research described in the article.

\medskip

%\section*{Acknowledgments}

\end{document}